\documentclass[11pt,a4paper]{article}
\usepackage{color}
\usepackage{graphicx}
\usepackage{amsfonts}
\usepackage{extarrows}
\usepackage{amsmath,amsthm,amssymb,color}
\usepackage{hyperref}
\usepackage{eepic}
\usepackage{lineno}
\usepackage{enumerate}	
\usepackage{paralist}
\usepackage{cite}
\usepackage{algorithm}
\usepackage{algorithmicx}
\usepackage{algpseudocode}

\usepackage{tikz}

\newtheorem{theorem}{Theorem}[section]
\newtheorem{proposition}[theorem]{Proposition}
\newtheorem{lemma}[theorem]{Lemma}

\theoremstyle{definition}

\newtheorem{claim}{Claim}

\newtheorem{conjecture}{Conjecture}[section]

\begin{document}
\title{\bf A strengthening on consecutive odd cycles in graphs of given minimum degree}
\author{Hao Lin \thanks{Department of Information Security, Naval University of Engineering, Wuhan, China. Supported by the Foundation of Naval University (2025508020).
Email: \texttt{haolinz6@qq.com}.}
\quad Guanghui Wang\thanks{ State Key Laboratory of Cryptography and Digital Economy Security, Shandong
University, Jinan, China, and School of Mathematics, Shandong University, 
Jinan, China. Supported by National Key R\&D Program of China (2020YFA0712400) and  Natural Science Foundation
of China (12231018). Email: \texttt{ghwang@sdu.edu.cn}.}
\quad Wenling Zhou\textsuperscript{\P}\thanks{School of Mathematics, Shandong University, Jinan,  China.  Supported by  Natural Science Foundation of China (12401457),  the China Postdoctoral Science Foundation (2024M761780)  and Natural Science Foundation of Shandong Province (ZR2024QA067).  Email: \texttt{gracezhou@sdu.edu.cn}. \\
\text{\textsuperscript{\P}Corresponding author.}}}

\date{}
\maketitle


\begin{center}
\begin{minipage}{130mm}
\small\noindent{\bf Abstract:}
Liu and Ma [{\it J. Combin. Theory Ser. B, 2018}] 
conjectured that every $2$-connected non-bipartite graph with minimum degree at least $k+1$ contains $\lceil k/2\rceil $ cycles with consecutive odd lengths. In particular, 
they showed that this conjecture holds when $k$ is even. In this paper, we confirm this conjecture for any $k\in \mathbb N$.
Moreover, we also improve some previous results about cycles of consecutive lengths.


\end{minipage}
\end{center}

\section{Introduction}

The research of the distribution of cycle lengths is a fundamental area in graph theory and
Erdős~\cite{Erdos-cycles-1976, Erdos-cycles-1992, Erdos-cycles-1995,Erdos-cycles-1997} posted many early problems on this topic. 
During the last decades, there has been extensive research on consecutive (even or odd) cycles in relation to minimum degree, average degree, connectivity, and chromatic number, see\cite{mim-degree-1984, Fan-2002, Liu-Ma2018,average-degree-2016,connectivity-2021,Gao-Li-2024, Gao-Ma-Siam-2021} (just to mention a few). 
In 2015, Liu and Ma~\cite{Liu-Ma2018} obtained several breakthrough results on the relation between cycle lengths and minimum degree. In particular, they proved the following general theorem:
for any positive integer $k$, every  $2$-connected non-bipartite graph $G$	with minimum degree at least $k+1$  contains $\lfloor k/2 \rfloor$ cycles with consecutive odd lengths, see {\cite[Theorem 1.3]{Liu-Ma2018}}.

Note that the above theorem is nearly tight, as the complete graph $K_{k+2}$  on $k+2$ vertices has exactly $\lceil k/2\rceil $ cycles with consecutive odd lengths. We quote from Liu and Ma~\cite{Liu-Ma2018} that ``It will be interesting if one can close the gap between our results and the best possible upper bounds." Therefore, they proposed the following conjecture.

\begin{conjecture}[{\hspace{-0.05em}\cite[Conjecture 6.1]{Liu-Ma2018}}] \label{Conj:Liu-Ma}
If $G$ is a $2$-connected non-bipartite graph with minimum degree at least $k+1$, then $G$ contains $\lceil k/2\rceil $ cycles with consecutive odd lengths. 
\end{conjecture}

In this paper, we confirm this conjecture and obtain the following.

\begin{theorem}\label{Thm:main}
Let $k\ge 1$ be an integer and $G$ be a $2$-connected non-bipartite graph. If the minimum degree of $G$  is at least $k+1$, then $G$ contains $\lceil k/2\rceil $ cycles with consecutive odd lengths. 
\end{theorem}

As mentioned above, the minimum degree condition in Theorem~\ref{Thm:main} is best possible. Moreover, Theorem~\ref{Thm:main} improves many previous results such as {\cite[Theorem 1.12]{Liu-Ma2018}}.

Another natural question that arises is: What are the necessary or sufficient conditions for the existence of $k$ cycles with consecutive lengths? Clearly, these conditions must involve non-bipartiteness when $k\ge 2$. Furthermore, Bondy and Vince~\cite{Bondy-1998} demonstrated that the 3-connectivity is necessary by constructing infinitely many $2$-connected non-bipartite graphs with arbitrarily large minimum degrees, yet not containing two cycles of consecutive lengths. In particular, they asked whether there exists a function $f(k)$ such that every $3$-connected non-bipartite graph with minimum degree at least $f(k)$ contains $k$ cycles with consecutive lengths. 
Fan~\cite{Fan-2002} proved the existence of $f(k)$ by proved $f(k)\le 3\lceil k/2\rceil $.
Liu and Ma~\cite{Liu-Ma2018} improved the result to $f(k)\le k+4$.
Recently, Gao, Huo, Liu and Ma~\cite{Gao-Ma-IMRN2022} showed that $f(k)\le k+1$. From the complete graph $K_{k+1}$, this minimum degree condition is tight. 
Our next result shows that, in fact, the minimum degree condition can be relaxed to $k$ except for $K_{k+1}$.

\begin{theorem}\label{Thm:k-consecutive-cycle}
Let $k\ge 6$ be an integer and $G$ be a $3$-connected non-bipartite graph. If the minimum degree of $G$ is at least $k$, then $G$ contains $k$ cycles of consecutive lengths, except that $G$ is $K_{k+1}$.
\end{theorem}

By considering the graph $K_{k+1}-M$ obtained from $K_{k+1}$ by deleting a non-empty matching $M$, we see that the minimum degree condition in Theorem~\ref{Thm:k-consecutive-cycle} is tight. 
In addition, applying Theorem~\ref{Thm:k-consecutive-cycle} with $k+1$, we can obtain $\lceil k/2\rceil $ cycles with consecutive odd lengths in $G$. 
Hence, we can present the proof of Theorem~\ref{Thm:main} by reducing to the following two results.

\begin{theorem}\label{Thm:case1}
Let $k\ge 1$ be an integer and $G$ be a $2$-connected non-bipartite graph. 
If $G$ is not $3$-connected and the minimum degree of $G$  is at least $k+1$, then $G$ contains $\lceil k/2\rceil $ cycles with consecutive odd lengths. 
\end{theorem}

\begin{theorem}\label{Thm:case2}
Every $3$-connected  non-bipartite graph $G$ with minimum degree at least four contains two cycles with consecutive odd lengths. 
\end{theorem}

Clearly, the case $k \in \{1, 2\}$ of Theorem~\ref{Thm:main} is trivial and the case $k \in \{3,4\}$ follows Theorem~\ref{Thm:case1} and Theorem~\ref{Thm:case2}. When $k\ge 5$,  the proof of Theorem~\ref{Thm:main} follows from Theorem~\ref{Thm:k-consecutive-cycle} and Theorem~\ref{Thm:case1}  immediately.

 The rest of this paper is organized as follows. 
In the next section, we will introduce some notations and results used in subsequent proofs.
In Section~\ref{sec:k-consecutive-cycle}, we will prove Theorem~\ref{Thm:k-consecutive-cycle} and Theorem~\ref{Thm:case2}, respectively.
In Section~\ref{sec:case1}, we will give a proof of Theorem~\ref{Thm:case1}.

\section{Notation and preliminaries}\label{sec-Preliminaries}

We remark that our notation follows~\cite{diestel2024graph}  and only considers simple undirected graphs. 
For a positive integer $k$, we denote by $[k]$ the set $\{1, \dots, k\}$.
Let $G$ be a graph and $S\subseteq V(G)$ be a vertex subset of $G$. 
We denote by $N_G (v)$ the neighborhood of a vertex $v$ in $G$.
We define $G[S]$ to be the subgraph induced by $S$ 
in $G$,  and $G-S$ to be the subgraph $G[V\setminus S]$.
For a vertex $v$ of $G$, the \emph{degree} of $v$, denoted by $d_G(v)$, is the number of edges in $G$ incident with $v$.
For two distinct vertices $x$, $y$ of $G$, let $G + xy$ 
denote the graph obtained from $G$ by adding the
edge $xy$, i.e., $E(G + xy)=E(G)\cup \{xy\}$;
and let $G-xy$ denote the graph obtained from $G$ by removing the edge $xy$, i.e., $E(G - xy)=E(G)\setminus  \{xy\}$.
A vertex subset $A$ of $G$ is called a \emph{cut-set} of $G$ if $G-A$ contains more components than $G$. 
If $A=\{v\}$ is a singleton, then we say that $v$ is a \emph{cut-vertex} of $G$.
A \emph{block} $B$ of $G$ is a maximal connected subgraph of $G$ such that $B$ contains no cut-vertex. So, a block is either an isolated vertex, an edge or a $2$-connected graph. 
An \emph{end-block} in $G$ is a block in $G$ containing at most one cut-vertex of $G$.
If $D$ is an end-block of $G$ and a vertex $x$ is the only cut-vertex of $G$ with $x \in V(D)$, then we say that $D$ is an \emph{end-block with cut-vertex $x$}.
For a cycle $C$, we use $\vec C$ to express this cycle with a prescribed orientation, and for vertices $u,v\in V(C)$, the notation $\vec{C}[u,v]$ denotes the subpath of $C$ from $u$ to $v$ following the given orientation. For a vertex set $\{v_1, v_2, \dots, v_t\}$, we use $(v_1, v_2, \dots, v_t)$ to express a path from $v_1$ to $v_t$ and passing through $v_2, \dots, v_{t-1}$ in turn.

We will appeal to a recent result of Gao, Huo, Liu and Ma~\cite{Gao-Ma-IMRN2022}, which
was used to resolve a number of conjectures regarding the existence of cycles of prescribed lengths.
We say that $k$ paths $P_1, P_2, \dots , P_k$ are \emph{admissible}, if $|E(P_1)| \ge 2$, and either $|E(P_{i+1})-E(P_i)|=1$ for all $i\in [k-1]$ or $|E(P_{i+1})-E(P_i)|=2$ for all $i\in [k-1]$. The following is a key tool in~\cite{Gao-Ma-IMRN2022}.

\begin{theorem}[{\hspace{-0.05em}\cite[Theorem 3.1]{Gao-Ma-IMRN2022}}]\label{Thm3.1:Gao-Ma-IMRN2022}
Let $k\ge 1$ be an integer and $G$ be a graph with $x, y \in V (G)$ such that $G+xy$ is $2$-connected. If every  $v\in V(G)\setminus\{x,y\}$ has degree at least $k +1$, then there exist $k$ admissible paths from $x$ to $y$ in $G-xy$.
\end{theorem}

Recently, Chiba, Ota and Yamashita~\cite{Chiba-2023} proved that the degree condition in Theorem~\ref{Thm3.1:Gao-Ma-IMRN2022} can be relaxed as follows.

\begin{theorem}[{\hspace{-0.05em}\cite[Theorem 3]{Chiba-2023}}]\label{thm:Chiba-2023-thm-1}
Let $k\ge 1$ be an integer, $G$ be a graph with $x, y \in V (G)$ such that $G+xy$ is $2$-connected, and $z$ be a vertex of $G$ {\rm (}possibly $z \in  \{x, y\}${\rm )}such that $V (G)\setminus\{x, y, z\} \neq \emptyset$. If every $v\in V(G)\setminus\{x,y,z\}$ has degree at least $k +1$, then there exist $k$ admissible paths from $x$ to $y$ in $G-xy$.
\end{theorem}

For the merging of consecutive paths and admissible paths, we have the following lemma.

\begin{lemma}\label{lem: merge-path}
	Let $s\ge 2$ and $t$ be two positive integers, and let $G$ be a graph with $x,y,z\in V(G)$. For a subset $W\subseteq V(G)\setminus \{z\}$ containing $x$ and $y$, if there exist $s$ consecutive paths $P_1,\dots, P_s$ from $x$ to $y$ in $G[W]$, and there exist $t$ admissible paths $Q_1, \dots, Q_t$ from $y$ to $z$ in $G-W\setminus\{y\}$. Then there exist at least $s+t-1$ consecutive paths from $x$ to $z$ in $G$.
\end{lemma}

\begin{proof}
Let $A$ be an arithmetic progression with common difference one.
	If $B$ is also an arithmetic progression with common difference one, then $A + B = \{a + b: a \in A, b \in B\}$ forms an arithmetic progression with common difference one of size $|A|+|B|-1$.
	If $B$ is an arithmetic progression with common difference two, then $A + B = \{a + b: a \in A, b \in B\}$ forms an arithmetic progression with common difference one of size $|A|+2(|B|-1)$.
	So the set $\{P_i\cup Q_j: i \in [s], j\in [t]\}$ contains at least $s + t-1$ consecutive paths from $x$ to $z$ in $G$.
\end{proof}

We also need a concept on cycles, which is crucial in the proof of Theorem~\ref{Thm:k-consecutive-cycle} and Theorem~\ref{Thm:case2}.
We say that a cycle $C$ in a connected graph $G$ is \emph{non-separating} if $G-V(C)$ is connected.
In $3$-connected graphs, since the non-separating induced cycles play a special role, the proof of the following lemma can be found in several groups (see~\cite{Bondy-1998,Thomassen-induced-cycles}), though it was not formally stated.

\begin{lemma}[{\hspace{-0.05em}\cite{Bondy-1998}}]\label{lem:Bondy-1998}
Every $3$-connected non-bipartite graph contains a non-separating induced odd cycle.
\end{lemma}

We also need the following lemma on non-separating odd cycles, which is implicitly proved in~\cite[Lemma 5.1]{Liu-Ma2018}. 
Thus, we omit its proof and refer the reader to~\cite[Lemma 5.1]{Liu-Ma2018}.

\begin{lemma}\label{lem:triangle}
Let $G$ be a graph with minimum degree at least four. If $G$ contains a non-separating induced odd cycle, then 
the shortest non-separating induced odd cycle $C$ in $G$,
denoted by $v_0v_1\dots v_{2s}v_0$,  satisfies either
\begin{itemize}
\item[{\rm ($\spadesuit$)}] $C$ is a triangle, or 

\item[{\rm ($\clubsuit$)}]  for every non-cut-vertex $v$ of $G-V(C)$, $|N_G(v) \cap V(C)| \le 2$, and the equality holds if and only if $N_G(v) \cap V(C) =\{v_i, v_{i+2}\}$ for some $i$, where the indices are taken under the additive group $\mathbb Z_{2s+1}$.
	\end{itemize}
\end{lemma}

\section{Proofs of Theorem~\ref{Thm:k-consecutive-cycle} and Theorem~\ref{Thm:case2} }\label{sec:k-consecutive-cycle}
This section is devoted to the proofs of Theorem~\ref{Thm:k-consecutive-cycle} and Theorem~\ref{Thm:case2}. 
Let $G$ be a $3$-connected non-bipartite graph with minimum degree at least $k$ with $k\ge 4$. By Lemma~\ref{lem:Bondy-1998} and Lemma~\ref{lem:triangle}, $G$ contains a non-separating induced odd cycle $C$ that satisfies the conclusion of Lemma~\ref{lem:triangle}.
Depending on whether $G$ contains a triangle, the proofs of Theorem~\ref{Thm:k-consecutive-cycle} and Theorem~\ref{Thm:case2} can be naturally divided into two cases. When  $G$ contains triangles,  Gao, Huo and Ma~\cite{Gao-Ma-Siam-2021} proved the following sharp result on consecutive cycles in graphs.

\begin{lemma}[{\hspace{-0.05em}\cite[Theorem 4.1]{Gao-Ma-Siam-2021}}]\label{Thm4.1:Gao-Ma-Siam-2021}
	Let $k \ge 2$ be an integer and $G$ be a $2$-connected graph containing a triangle. If the minimum degree of $G$ is at least $k$, then $G$  contains $k$ cycles of consecutive lengths, except that $G = K_{k+1}$.
\end{lemma}

Based on the above result, it suffices to consider the case where $G$ is triangle-free in the subsequent proofs of Theorem~\ref{Thm:k-consecutive-cycle} and Theorem~\ref{Thm:case2}. Combined with Lemma~\ref{lem:triangle}, we can further assume that the shortest non-separating induced odd cycle $C$ in $G$, denoted by $v_0v_1 \dots v_{2s}v_0$ for some $s \geq 2$, satisfies the property ($\clubsuit$).

Throughout this section, the subscripts will be taken under the additive group $\mathbb{Z}_{2s+1}$. Let $G_1:=G-V(C)$. 
We first have the following claim.

\begin{claim}\label{claim:degree=k-1}
Given $k \ge 4$,  let $B =G_1$ if  $G_1$ is $2$‐connected; otherwise, let $B$ be an  end-block of $G_1$ with cut-vertex $x$. If there exists $v\in V(B)\setminus\{x\}$ such that $|N_G(v) \cap V(C)|=2$, then $G$  contains $k$ cycles of consecutive lengths.
\end{claim}
\begin{proof}
Note that every vertex in $V(B)\setminus\{x\}$ is not a cut-vertex of $G_1$. By property ($\clubsuit$), we see that every vertex of $ V(B)\setminus\{x\}$ has degree at least $k-2$ in $B$. 
Without loss of generality, suppose that there exists $v\in V(B)\setminus\{x\}$ such that $N_G(v) \cap V(C)=\{v_{2s}, v_{1}\}$.
Since the minimum degree of $G$ is at least $k\ge 4$, we can choose a vertex $u\in N_G(v_s)\cap V(G_1)$ such $u\neq v$. Note that $(v, v_1, v_2, \dots, v_s, u)$, $(v, v_{2s}, v_{2s-1}, \dots, v_s, u)$, $(v, v_{2s}, v_0, v_{1}, v_{2}, \dots, v_s, u)$ and $(v, v_1, v_0, v_{2s}, v_{2s-1}, \dots, v_s, u)$ are $4$ consecutive paths from $v$ to $u$ of lengths $s+1$, $s+2$, $s+3$ and $s+4$, respectively. If $u\in V(B)$, 
then applying Theorem~\ref{thm:Chiba-2023-thm-1} to $B$ with $u,v,x\in V(B)$, there exist $k-3$ admissible paths from $v$ to $u$ in $B$. If $u\notin V(B)$,  then we can find a path $T_{ux}$ from $u$ to $x$ such that $E(T_{ux})\cap E(B)=\emptyset$ since $G_1$ is connected. Then, we apply Theorem~\ref{Thm3.1:Gao-Ma-IMRN2022} to $B$ with $x,v \in V(B)$ and find that $k-3$ admissible paths from $x$ to $v$ in $B$. Therefore, we can always find $k-3$ admissible paths from $v$ to $u$ in $G_1$. 
By Lemma~\ref{lem: merge-path}, merging these $4$ consecutive paths and $k-3$ admissible paths produces at least $4+(k-3)-1=k$ cycles in $G$ with consecutive lengths.  
\end{proof}

Based on the above discussions, we prove Theorem~\ref{Thm:k-consecutive-cycle} and Theorem~\ref{Thm:case2} in the following subsections, respectively.
\subsection{Proof of Theorem~\ref{Thm:k-consecutive-cycle}}

\begin{proof}[{\bf Proof of Theorem~\ref{Thm:k-consecutive-cycle}}]
Let  $k\ge 6$ be an integer, $G$ be a $3$-connected non-bipartite graph with minimum degree at least $k$, and every shortest non-separating induced odd cycle of $G$ satisfies the property ($\clubsuit$). 
Depending on whether  $G_1$ is $2$-connected,  we divide the rest of the proof into two cases.
\smallskip

{\bf Case 1. } $G_1$  is $2$-connected.
\smallskip

By Claim~\ref{claim:degree=k-1}, 
 we can assume that every vertex $v\in V(G_1)$ satisfies that $|N_G(v) \cap V(C)| \le 1$ which imples that the minimum degree of $G_1$ is at least $k-1$. We choose a vertex $v\in N_G(v_0)\cap V(G_1)$ and a vertex $u\in N_G(v_{s-1})\cap V(G_1)$ such that $u\neq v$. Note that $(v, v_0, v_1,v_2, \dots, v_{s-1}, u)$ and $(v, v_{0}, v_{2s}, v_{2s-1}, \dots, v_{s-1}, u)$ are two paths from $v$ to $u$ of lengths $s+1$ and $s+4$, denoted $Q_1$  and $Q_2$, respectively. 
 Applying Theorem~\ref{Thm3.1:Gao-Ma-IMRN2022} to $G_1$ with $u,v\in V(G_1)$, then $G_1$ contains $k-2$ admissible paths $P_1, P_2, \dots, P_{k-2}$ from $v$ to $u$. If $|E(P_{j+1})-E(P_{j})|=1$ for all $j\in [k-3]$, then the set $\{Q_i\cup P_j: i\in [2], j\in [k-2]\}$ contains  $k+1$ consecutive cycles in $G$ when $k\ge 5$.
If $|E(P_{j+1})-E(P_{j})|=2$ for all  $j\in [k-3]$,  then the set $\{Q_i\cup P_j: i\in [2], j\in [k-2]\}$ contains  $2(k-3)\ge k$ consecutive cycles in $G$ when $k\ge 6$. 
\smallskip

{\bf Case 2. } $G_1$ is not $2$-connected.
\smallskip

 By Claim~\ref{claim:degree=k-1}, 
 we can assume that for any end-block $B$ of $G_1$, every vertex of $B$ other than the cut-vertex has degree at least $k-1$.
Let $D_1$ be an end-block with cut-vertex $x$ of $G_1$ and $G_2:=V(G_1)-(V(D_1)\setminus \{x\})$. 
Note that both $D_1$ and $G_2$ are neither an isolated vertex nor an edge.
We first have the following claim.

\begin{claim}\label{claim:quasi‐diagonal}
There exists an $i\in \mathbb Z_{2s+1}$ such that $N_G(v_i)\cap (V(D_1)\setminus \{x\})\neq \emptyset $ and $N_G(v_{i+s})\cap V(G_2)\neq \emptyset $. 
\end{claim}
\noindent{ \emph{Proof}.} 
For simplicity, let $A:=V(D_1)\setminus \{x\}$ and $B:=V(G_2)$. 
For the sake of contradiction, assume that each $i\in \mathbb Z_{2s+1}$ satisfies that either $N_G(v_i)\cap A = \emptyset $ or $N_G(v_{i+s})\cap B= \emptyset $. Fix $i\in \mathbb{Z}_{2s+1}$. If $N_G(v_i)\cap A=\emptyset$, then we must have $N_G(v_i)\cap B\neq \emptyset$. By the assumption, we further have  $N_G(v_{i+s})\cap A=\emptyset$.
We now construct an auxiliary graph $H_C$ on $V (C)$ with $E(H_C)=\{v_iv_{i+s} : i\in \mathbb{Z}_{2s+1}\}$. Then $H_C$ is isomorphic to $C$ since $s$ and $2s+1$ are relatively prime.  By the cyclicity of $H_C$, $N_G(v_{i})\cap A=\emptyset$ for all $i\in \mathbb Z_{2s+1}$ , then $x$ is a cut-vertex of $G$, 
a contradiction.
\hfill $\blacksquare$
\smallskip

Now let $i' \in \mathbb Z_{2s+1}$ such that $N_G(v_{i'})\cap (V(D_1)\setminus \{x\})=\{v\}$ and $N_G(v_{i'+s})\cap V(G_2)=\{u\}$.  
Note that $D_1$ is $2$-connected and every vertex of $D_1$ other than $x$ has degree at least $k-1$. By Theorem~\ref{Thm3.1:Gao-Ma-IMRN2022}, $D_1$ contains $k-2$ admissible paths $P_1, \dots, P_{k-2}$ from $x$ to $v$ in $D_1$.  
Depending on the distribution of lengths of $P_1, \dots, P_{k-2}$, we divide the remainder of the proof into two subcases.
\smallskip

{\bf Case 2.1. } $|E(P_{j+1})-E(P_{j})|=2$ for all $j\in [k-3]$.
\smallskip

Note that $\vec{C}[v_{i'},v_{i'+s}]$ and $\vec{C}[v_{i'+s},v_{i'}]$ are two consecutive subpath of $C$. Thus, we can get two consecutive paths $Q_1, Q_2$ from $v$ to $u$ such that all internal vertices of $Q_1$ and $Q_2$ belong to $V(C)$.
Let $T$ be a path from $x$ to $u$ in $G_2$, and possibly $x=u$.
In this case, the set $\{Q_i\cup T \cup P_j: i\in [2], j\in [k-2]\}$ contains $2(k-2)\ge k$ consecutive cycles in $G$ when $k\ge 4$.
\smallskip

{\bf Case 2.2. } $|E(P_{j+1})-E(P_{j})|=1$ for all $j\in [k-3]$.
\smallskip

 In this scenario, let $D_2\subseteq G_2$ be an end-block of $G_1$ with cut-vertex $y$, and possibly $x=y$.  we first choose a path $P_{xy}$ from $x$ to $y$ such that $E(P_{xy})\cap (E(D_1) \cup E(D_2))=\emptyset$.
Since $D_2$ is an end-block of $G_1$ and $G$ is $2$-connected, there exists a vertex $z\in V(D_2)\setminus\{y\}$ such that there exists $w\in N_G(z)\cap V(C)$. 
Hence, there is a path $P_{vz}$ from $v$ to $z$ such that $V(P_{vz})\subset V(C)\cup \{v,z\}$. 
Secondly, $D_2$ is $2$-connected and every vertex of $D_2$ other than $y$ has degree at least $k-1$. 
 By Theorem~\ref{Thm3.1:Gao-Ma-IMRN2022}, $D_2$ contains $k-2$ admissible paths $P'_1, \dots,  P'_{k-2}$ from $z$ to $y$ in $D_2$. 
Recalling that $P_1, \dots, P_{k-2}$ are $k-2$ consecutive paths from $x$ to $v$ in $D_1$.
So, the set $\{P_j\cup P_{vz}\cup P'_{\ell}\cup P_{xy}: j\in [k-2], \ell \in [k-2]\}$ contains $2(k-2)-1\ge k$ consecutive cycles in $G$ when $k\ge 5$.
This completes the proof of Theorem~\ref{Thm:k-consecutive-cycle}.
\end{proof}

We remark that the above proof works for the following result as well,  which is a slight extension of~\cite[Theorem 4.4]{Gao-Ma-IMRN2022}.

\begin{theorem}
Let $k\ge 6$ be an integer and $G$ be a $2$-connected graph containing a non-separating induced odd cycle. If the minimum degree of $G$ is at least $k$, then $G$ contains $k$ cycles of consecutive lengths, except that $G$ is $K_{k+1}$.
\end{theorem}

\subsection{Proof of Theorem~\ref{Thm:case2}}\label{sec:case2}
In this subsection, we provide the proof of Theorem~\ref{Thm:case2}. We begin by introducing a useful result.

\begin{lemma}[{\hspace{-0.05em}\cite[Theorem 5]{Gao-Li-2024}}]\label{lem:Gao-Li-thm-5}
Let $G$ be a graph with $x, y \in V (G)$ such that $G + xy$ is $2$‐connected. If every $v\in V(G)\setminus\{x,y\}$ has degree at least $3$, and any edge $uv \in E (G)$ with $\{u, v\} \cap \{x, y\} = \emptyset$ has degree sum $d_G (u)+d_G (v) \ge  7$, then there exist two paths from $x$ to $y$ in $G  -xy$ whose lengths differ by two.
\end{lemma}

\begin{proof}[{\bf Proof of Theorem~\ref{Thm:case2}}]
Let $G$ be a $3$-connected  non-bipartite graph with minimum degree at least four, and every shortest non-separating induced odd cycle in $G$ satisfies the property ($\clubsuit$). 
 We aim to show that $G$ contains two cycles with consecutive odd lengths. Suppose that $G$ is a counterexample.
 Depending on whether  $G_1$ is $2$-connected,  we divide the rest of the proof into two cases.
\smallskip

{\bf Case 1. } $G_1$ is not $2$-connected.
\smallskip

 Let $D_1$ and $D_2$ are two distinct end-blocks of $G_1$ with cut-vertices $x_1$ and $x_2$ (possibly $x_1=x_2$), respectively. 
 Since $G$ is $3$-connected, there exist two distinct vertices $v_p, v_q\in V(C)$ such that there is $u_1\in N_G(v_p)\cap (V(D_1)\setminus\{x_1\})$ and $u_2\in N_G(v_q)\cap (V(D_2)\setminus\{x_2\})$. By Claim~\ref{claim:degree=k-1}, 
 every vertex in $V(D_i)\setminus
\{x_i\}$ has degree at least three in $D_i$ for all $i\in [2]$. By Theorem~\ref{Thm3.1:Gao-Ma-IMRN2022}, $D_1$ contains two admissible paths $P_1$ and $P_2$ from $x_1$ to $u_1$ and $D_2$ contains two admissible paths $P'_1$ and $P'_2$ from $x_2$ to $u_2$. 
Note that $u_1$ and $u_2$ are cut-vertices of $G_1$, we can choose a path $T_x$ from $x_1$ to $x_2$ in $G_1$ such that $E(T_{x})\cap (E(D_1)\cup E(D_2))=\emptyset$. 
Thus, by concatenating $P_i$, $T_x$ and $P'_j$ for all $i,j\in [2]$, we can obtain at least three admissible paths $P^*_1, P^*_2, P^*_3$ from $u_1$ to $u_2$ in $G_1$.  
Note that the lengths of $\vec{C}[v_{p},v_{q}]$ and $\vec{C}[v_{q},v_{p}]$  have different parity.
Hence, the family of $\vec{C}[v_{p},v_{q}]\cup P^*_\ell$ and $\vec{C}[v_{q},v_{p}]\cup P^*_\ell$ for all $\ell\in [3]$ must contain 
three cycles whose lengths differ by one or two,  and the shortest one of which is odd, a contradiction.
\smallskip

{\bf Case 2. } $G_1$ is $2$-connected.
\smallskip

In this case, we first have the following claim.
\begin{claim}\label{cliam:degree-sum}
There is an edge $uv\in E(G_1)$ such that both $u$ and $v$ are adjacent to the cycle $C$. 
\end{claim}
\noindent {\emph{Proof}. }
By Claim~\ref{claim:degree=k-1}, we first have $d_C(u)\le 1$ for each $u\in V(G_1)$.
Suppose in contrast that $d_C(u)+d_C(v)\le 1$ for any $uv\in E(G_1)$. Then by the minimum degree of $G$, $d_{G_1}(u)+d_{G_1}(v)\ge 7$. In addition, we can choose two  vertex-disjoint edges $x'x, y'y\in E(G)$ such that $x',y'\in V(C)$ and $x, y\in V(G_1)$. In this case, we can apply Lemma~\ref{lem:Gao-Li-thm-5} to $G_1$ with $x,y\in V(G_1)$, and obtain two paths from $x$ to $y$ in $G_1$ whose lengths differ by two.
Since $\vec{C}[x', y']$ and $\vec{C}[y',x']$ have lengths with different parity, we can find two cycles of consecutive odd lengths in $G$, a contradiction.
\hfill $\blacksquare$
\smallskip

By Cliam~\ref{cliam:degree-sum}, we define 
$
\mathcal {F}=\{\{a,b\}\subset V(G_1): ab\in E(G_1), |N_G(a) \cap V(C)|=|N_G(b) \cap V(C)|= 1\}
$.
For any $\{a,b\}\in \mathcal {F}$, let $\{u_a\}:= N_G(a)\cap V(C)$ and $\{u_b\}:= N_G(b)\cap V(C)$. Note that $u_a\neq u_b$ since $G$ is triangle-free.
Without loss of generality, let $\vec C[u_a, u_b]$ be the even path and $\vec C[u_b, u_a]$ be the odd path on $C$. 
Set $C_{ab}:=\vec C[u_a, u_b]\cup \{u_aa, ab, bu_b\}$. By the above definitions, we have the following observation: For any $\{a,b\}\in \mathcal {F}$, $C_{ab}$ is an induced odd cycle of $G$. Otherwise,  $u_au_b\in E(G)$ and $C$ and $C_{ab}$ are two cycles of consecutive odd lengths in $G$, a contradiction. 

Next, depending on the property of $\{a,b\}$  in $G_1$, we divide the rest of the proof into two subcases.
\smallskip

{\bf Case 2.1.} There exists $\{a,b\}\in \mathcal {F}$ that is not a cut-set of $G_1$. 
\smallskip

In this case, the induced odd cycle $C_{ab}$ is non-separating in $G$. By the minimality of $C$, the odd path $\vec C[u_b, u_a]$ only has four vertices, denoted by $\vec C[u_b, u_a]=(u_b, u^{-}_b, u^+_a, u_a)$. 
Clearly, $C_{ab}$ is also a shortest non-separating induced odd cycle of $G$, which implies that $G-V(C_{ab})$ is also $2$-connected with minimum degree at least three; otherwise,  by replacing $C$ with $C_{ab}$ in the proof of {\bf Case 1}, we can find two cycles of consecutive odd lengths in
$G$, a contradiction.

If $G_1-\{a,b\}$ is $2$-connected, let $B = G_1-\{a,b\}$; otherwise, let $B$ be an end-block  of $G_1-\{a,b\}$ with cut-vertex $x$. Note that $B$ is not a single vertex since $G$ is triangle-free.  Moreover, we see that $(N_G(a)\cup N_G(b))\cap (V(B)\setminus \{x\})\neq \emptyset$ since $G_1$ is $2$‐connected.
Without loss of generality, suppose that there exists $v_a\in V(B)\setminus \{x\}$ such that $av_a\in E(G)$. Let $G_2:=G[V(B)\cup \{u^+_a, a\}]$. Then we have the following claim.

\begin{claim}\label{cliam:odd-cycle}
For any $v\in V(B)\setminus\{x\}$, $d_{G_2}(v)\ge 3$.
\end{claim}
\noindent {\emph{Proof}. } 
Note that for any $v\in V(B)\setminus\{x\}$, $N_G(v)\subset V(G_2)\cup \{u^-_b, b\}\cup V(\vec C[u_a, u_b])$ and $d_G(v)\ge 4$.
If there exists  $v^*\in V(B)\setminus\{x\}$ such that the set $N_G(v^*)\cap (V(G)\setminus V(G_2))$ contains at least two vertices. Then the set $N_G(v^*)\cap (V(G)\setminus V(G_2))$  is exactly $\{u^-_b, b\}$ since $|N_G(v^*) \cap V(C_{ab})| \le 1$ and $|N_G(v^*) \cap V(C)| \le 1$.   In this case, $u_bb$ and $P:=(u_b,u^-_b,v^*,b)$ form two paths from $u_b$ to $b$ of consecutive odd lengths. Thus, by concatenating the path $C_{ab}- u_bb$, we can obtain two cycles of consecutive odd lengths in $G$, a contradiction. 
\hfill $\blacksquare$
\smallskip

By Claim~\ref{cliam:odd-cycle}, we have that $d_{G_2}(v)\ge 3$ for any $v\in V(B)\setminus\{x\}$. 
In this case, if there exists one vertex $v^+_a$ in $V(B)\setminus \{ v_a\}$ such that $u^+_av^+_a\in E(G)$. Then
$G_2+u^+_aa$ is $2$-connected.
Applying Theorem~\ref{thm:Chiba-2023-thm-1} to $G_2$  with $\{u^+_a, a, x\}$, there are two admissible paths $P_1$ and $P_2$ from $u^+_a$ to $a$ in $G_2$. Let $Q_1:=(u^+_a, u_a, a)$, $Q_2:=(u^+_a,u^-_b, u_b,b, a)$ and $Q_3:=u^+_au_a\cup \vec C[u_a, u_b]\cup u_bb\cup ba$. Clearly, $Q_1$ and $Q_2$ are two consecutive even paths, and $Q_3$ is an odd path. 
By merging these paths $P_i$ and $Q_j$ for all $i\in [2]$ and $j\in [3]$,  we can find two cycles of consecutive odd lengths in $G$, a contradiction.

Therefore, $(V(B)\setminus \{v_a\})\cap N_G(u^+_a)=\emptyset$. In addition,  $v_a\notin N_G(u^+_a)$. Otherwise, $u^+_au_a$ and $P':=(u^+_a, v_a, a, u_a)$ form two paths from $u^+_a$ to $u_a$ of consecutive odd lengths. Thus, by concatenating the path $C- u^+_au_a$, we can obtain two cycles of consecutive odd lengths in $G$, a contradiction. 
In addition, $a\notin N_G(u^+_a)$ since $G$ is triangle-free. By $d_G(u^+_a)\ge 4$, 
we can find a path $P_{u^+_ax}$ from $u^+_a$ to $x$ in $G_1-\{a,b\}$ such that $E(P_{u^+_ax})\cap E(B)=\emptyset$. Let $G'_2:=G[V(B)\cup \{ a\}]$. 
We apply Theorem~\ref{Thm3.1:Gao-Ma-IMRN2022} to $G'_2$ with $x, a\in V(G'_2)$, then there are two admissible paths $P'_1$ and $P'_2$ from $x$ to $a$ in $G'_2$. 
Similarly, we can also find two cycles of consecutive odd lengths in $G$, a contradiction. 
\smallskip

{\bf Case 2.2.} For any $\{a,b\}\in \mathcal {F}$,  $\{a,b\}$ is a cut-set of $G_1$.  
\smallskip

We take a pair $\{a_1,b_1\}\in \mathcal {F}$, and let $D_1$ and $D_2$ be two component of $G_1-\{a_1, b_1\}$.
Note that for $i\in[2]$, $D_i$ has no isolated vertices since $G$ is triangle-free.
Set $H_1:=G[V(D_1)\cup \{a_1, b_1\}]$. 
 Moreover, we have the following claim. \begin{claim}\label{cliam:degree-sum-2}
There exists $a'b'\in E(D_1)$ such that $\{a', b'\} \in \mathcal {F}$.
\end{claim}
\noindent {\emph{Proof}.} 
Suppose that any $uv\in E(D_1)$ such that $\{a', b'\} \notin \mathcal {F}$, which means that  $
|N_G(u) \cap V(C)|+|N_G(v) \cap V(C)|\le 1$ by Cliam~\ref{claim:degree=k-1}. Note that $H_1$ is $2$-connected since $G_1$ is $2$-connected. In addition, every $v\in V(H_1)\setminus\{a_1, b_1\}$ has degree at least $3$, and any edge $uv \in E (H_1)$ with $\{u, v\} \cap \{a_1, b_1\} = \emptyset$ has degree sum $d_{H_1} (u)+d_{H_1} (v) \ge  7$. By Lemma~\ref{lem:Gao-Li-thm-5},  there are two paths from $a_1$ to $b_1$ in ${H_1} $ whose lengths differ by two. Therefore, similar to the proof of Claim~\ref{cliam:degree-sum},
 we can find two cycles of consecutive odd lengths in $G$, a contradiction.
\hfill $\blacksquare$
\smallskip

By Claim~\ref{cliam:degree-sum-2}, we take $a'b'\in E(D_1)$ such that $\{a', b'\} \in \mathcal {F}$. 
In this case, $\{a', b'\}$ is a cut-set of $G_1$. Recalling that
$\{u_{a_1}\}= N_G(a_1)\cap V(C)$, $\{u_{b_1}\}= N_G(b_2)\cap V(C)$, 
$\{u_{a'}\}= N_G(a')\cap V(C)$ and $\{u_{b'}\}= N_G(b')\cap V(C)$. Since $G$ is triangle-free, we have $u_{a_1}\neq u_{b_1}$ and $u_{a'}\neq u_{a'}$. Thus, we can assume that $u_{a_1}\neq u_{a'}$. Then we can find an even path $Q_e$ and an odd path $Q_o$ from ${a_1}$ to $a'$ with all internal vertices in $C$. 

Next, let $D'_1$ be the component of $H_1-\{a', b'\}$ that does not contain $\{a_1, b_1\}$, $H'_1:=G[V(D'_1)\cup \{a', b'\}]$, $H_2:=G[V(D_2)\cup \{a_1, b_1\}]$ and $H_3:=H_1-V(D'_1)$. 
Clearly, $H'_1, H_2, H_3$ are $2$-connected, and $\{a_1, b_1, a', b'\}\subseteq V(H_3)$. Thus, we can apply  Theorem~\ref{Thm3.1:Gao-Ma-IMRN2022} to $H_2$ with $a_1, b_1\in V(H_2)$, and obtain two admissible paths $P_1$ and $P_2$ from $a_1$ to $b_1$ in $H_2$. Similarly, we can find two admissible paths $P'_1$ and $P'_2$ from $a'$ to $b'$ in $H'_1$. 
Moreover, applying Menger’s theorem to $H_3$ with $\{a_1, b_1\}$ and $\{a', b'\}$, we can find two internal disjoint  paths from $\{a_1, b_1\}$ to  $\{a', b'\}$ in $H_3$. Thus, we can assume that there is a path $T_b$ from $b'$ to  $b_1$ in $H_3$ such that $V(T_b)\cap \{a_1, a'\}=\emptyset$.
Now note that the set  $\{Q_e\cup P'_i\cup T_b\cup P_j\}$ and $\{Q_o\cup P'_i\cup T_b\cup P_j\}$  for all $i, j\in [2]$ must contain 
three cycles $C_1, C_2, C_3$ whose lengths differ by one or two,  and the shortest one of which is odd, a contradiction. This completes the proof of Theorem~\ref{Thm:case2}.
\end{proof}

\section{Proof of Theorem~\ref{Thm:case1}}\label{sec:case1}
In this section, we conclude with the proof of Theorem~\ref{Thm:case1}. We begin with a succinct proposition.

\begin{proposition}\label{pro:odd-even}
	Let  $G$ be a non-bipartite graph. For any distinct vertices $x$ and $y$ of $G$, if $G+xy$ is $2$-connected, then there exists an odd path and an even path between $x$ and $y$ in $G$.
\end{proposition}

\begin{proof}
	Since $G$ is non-bipartite, there is an odd cycle $C$ in $G$.
	For any distinct vertices $x$ and $y$ of $G$, since $G+ xy$ is $2$-connected, then there is at least one path from $x$ to $y$ in $G-xy$. Let $P$
	be a such path from $x$ to $y$ such that the size of $V(P)\cap V(C)$ is maximum.
	
	Since $G+ xy$ is $2$-connected, we have $|V(P)\cap V(C)|\ge 2$. Scan $P$ from $x$ to $y$, let $a$ be the first vertex in $C$, and $b$ be the last vertex in $C$.  Clearly, there are two distinct paths between $a$ and $b$ in $C$. Besides, the length of one is odd and the length of the other is even. Therefore, we conclude that there always exists an odd path and an even path between $x$ and $y$ in $G$.
\end{proof}

We are ready to prove Theorem~\ref{Thm:case1} using Proposition~\ref{pro:odd-even} and Theorem~\ref{Thm3.1:Gao-Ma-IMRN2022}.

\begin{proof}[ {\bf Proof of Theorem~\ref{Thm:case1}}]
Since $G$ is $2$-connected but not $3$-connected, there exist two distinct vertices $x,y\in V(G)$ such that $G-\{x, y\}$ contains at least two components.
	Let $C$ be an arbitrary component of $G-\{x, y\}$. Define $G_1:=G[V(C)\cup \{x,y\}]$ and $G_2:=G[V(G) \setminus V(C)]$.
	Note that for any $i\in [2]$, $G_i+xy$ is $2$-connected and every vertex of $G_i$ other than $x$ and $y$ has degree at least $k + 1$.
	Now we distinguish between two cases.
\smallskip

	{\bf Case 1. } At least one of $G_1$ and $G_2$ is non-bipartite.
	\smallskip
	
	Without loss of generality we assume that $G_1$ is non-bipartite.  By Proposition~\ref{pro:odd-even}, there exists an odd path $P_o$ and an even path $P_e$ between $x$ and $y$ in $G_1$. Then we apply Theorem~\ref{Thm3.1:Gao-Ma-IMRN2022} to $G_2$ with $x,y\in V(G_2)$ and obtain $k$ admissible paths $P_1, \dots, P_k$ from $x$ to $y$ in $G_2-xy$. If $|E(P_1)|$ is odd, then the set $\{P_i\cup P_e: i\in [k]\}$ contains at least $\lceil k/2\rceil $ cycles with consecutive odd lengths. If $|E(P_1)|$ is even, then the set $\{P_i\cup P_o: i\in [k]\}$ contains at least $\lceil k/2\rceil $ cycles with consecutive odd lengths.
	\smallskip
	
	{\bf Case 2. } Both $G_1$ and $G_2$ are bipartite.
\smallskip

 In this case, we claim that either $x, y$ are in the same part of $G_1$ and $x, y$ are in distinct parts of $G_2$, or $x, y$ are in the same part of $G_2$ and $x, y$ are in distinct parts of $G_1$. Suppose to the contrary that $x, y$ are in the same part of $G_i$ for $i\in [2]$.   Let $A$ be the union of the part containing $x, y$ in $G_1$ and the part containing $x, y$ in $G_2$.  Let $B=V(G)\setminus A$. Since $C$ be a component of $G-\{x, y\}$, $\{A, B\}$ forms a bipartition of $G$, a contradiction. 
	If $x, y$ are in distinct parts of $G_i$ for $i\in [2]$. 
	Let $A'$ be the union of the part containing $x$ in $G_1$ and the part containing $x$ in $G_2$. Let $B'=V(G)\setminus A'$. Then $\{A', B'\}$  also forms a bipartition of $G$, a contradiction. This proves the claim. 
	
	Without loss of generality we assume that $x, y$ are in the same part of $G_1$ and $x, y$ are in distinct parts of $G_2$. Then every path from $x$ to $y$ in $G_1$ is even and every path from $x$ to $y$ in $G_2$ is odd. Let $P$ be a path from $x$ to $y$ in $G_1$.
	Then by Theorem~\ref{Thm3.1:Gao-Ma-IMRN2022}, there exist $k$ paths from $x$ to $y$ with consecutive odd lengths in $G_2-xy$. 
	Concatenating each of these paths with $P$, we obtain $k$ cycles of consecutive odd lengths.
	This finishes the proof of Theorem~\ref{Thm:case1}.
\end{proof}

\section{Concluding Remarks}

In this paper,  we obtained some new results on the relation between cycle lengths and minimum degree.
In particular, Theorem~\ref{Thm:k-consecutive-cycle} states that for $k\ge 6$, every $3$-connected non-bipartite graph $G$ with minimum degree at least $k$ contains $k$ cycles of consecutive lengths, except that $G$ is $K_{k+1}$. 
Actually, the proof of Theorem~\ref{Thm:k-consecutive-cycle} works for some special cases when $k=4$ or $k=5$. 
Thus, we conjecture that the condition 
$k\ge 6$ may be relaxed.

\begin{conjecture}
Let $k\ge 4$ be an integer. If $G$ is a $3$-connected non-bipartite graph with minimum degree at least $k$, then $G$ contains $k$ cycles of consecutive lengths, except when $G$ is $K_{k+1}$.
\end{conjecture}

We remark that the bound $k \ge 4$ in the above conjecture cannot be further relaxed. For instance, when $k = 3$, the Petersen graph serves as a counterexample. It is a $3$-connected non-bipartite graph with minimum degree at least $3$, yet it contains cycles only of lengths $5$, $6$, $8$, and $9$, and hence does not contain $3$ cycles of consecutive lengths.

\section*{Acknowledgements}
Sincere thanks to Jie Ma for his helpful discussions and for carefully reading a draft. We also acknowledge the Extremal Combinatorics Workshop held at Shandong University, Weihai, in August 2024, which inspired part of this research. In addition, the authors thank Yandong Bai for his helpful discussions, and the anonymous referee for carefully reading the manuscript and providing many helpful comments.

\bibliographystyle{abbrv}
\bibliography{ref}

\begin{thebibliography}{10}

\bibitem{Bondy-1998}
J.~A. Bondy and A.~Vince.
\newblock Cycles in a graph whose lengths differ by one or two.
\newblock {\em J. Graph Theory}, 27(1):11--15, 1998.

\bibitem{Chiba-2023}
S.~Chiba, K.~Ota, and T.~Yamashita.
\newblock Minimum degree conditions for the existence of a sequence of cycles
  whose lengths differ by one or two.
\newblock {\em J. Graph Theory}, 103(2):340--358, 2023.

\bibitem{diestel2024graph}
R.~Diestel.
\newblock {\em Graph theory}.
\newblock Graduate Texts in Mathematics. Springer, Berlin, fifth edition, 2018.

\bibitem{Erdos-cycles-1976}
P.~Erd\H{o}s.
\newblock Some recent problems and results in graph theory, combinatorics and
  number theory.
\newblock In {\em Proceedings of the {S}eventh {S}outheastern {C}onference on
  {C}ombinatorics, {G}raph {T}heory, and {C}omputing ({L}ouisiana {S}tate
  {U}niv., {B}aton {R}ouge, {L}a., 1976)}, volume No. XVII of {\em Congress.
  Numer.}, pages 3--14. Utilitas Math., Winnipeg, MB, 1976.

\bibitem{Erdos-cycles-1992}
P.~Erd\H{o}s.
\newblock On some of my favourite problems in various branches of
  combinatorics.
\newblock In {\em Fourth {C}zechoslovakian {S}ymposium on {C}ombinatorics,
  {G}raphs and {C}omplexity ({P}rachatice, 1990)}, volume~51 of {\em Ann.
  Discrete Math.}, pages 69--79. North-Holland, Amsterdam, 1992.

\bibitem{Erdos-cycles-1995}
P.~Erd\H{o}s.
\newblock Some of my favourite problems in number theory, combinatorics, and
  geometry.
\newblock In {\em Combinatorics Week (Portuguese) (S\~ao Paulo, 1994)},
  volume~2, pages 165--186, 1995.

\bibitem{Erdos-cycles-1997}
P.~Erd\H{o}s.
\newblock Some old and new problems in various branches of combinatorics.
\newblock In {\em Graphs and combinatorics (Marseille, 1995)}, volume 165/166,
  pages 227--231, 1997.

\bibitem{Fan-2002}
G.~Fan.
\newblock Distribution of cycle lengths in graphs.
\newblock {\em J. Combin. Theory Ser. B}, 84(2):187--202, 2002.

\bibitem{Gao-Ma-IMRN2022}
J.~Gao, Q.~Huo, C.-H. Liu, and J.~Ma.
\newblock A unified proof of conjectures on cycle lengths in graphs.
\newblock {\em Int. Math. Res. Not. IMRN}, (10):7615--7653, 2022.

\bibitem{Gao-Ma-Siam-2021}
J.~Gao, Q.~Huo, and J.~Ma.
\newblock A strengthening on odd cycles in graphs of given chromatic number.
\newblock {\em SIAM J. Discrete Math.}, 35(4):2317--2327, 2021.

\bibitem{Gao-Li-2024}
J.~Gao, B.~Li, J.~Ma, and T.~Xie.
\newblock On two cycles of consecutive even lengths.
\newblock {\em J. Graph Theory}, 106(2):225--238, 2024.

\bibitem{mim-degree-1984}
A.~Gy\'arf\'as, J.~Koml\'os, and E.~Szemer\'edi.
\newblock On the distribution of cycle lengths in graphs.
\newblock {\em J. Graph Theory}, 8(4):441--462, 1984.

\bibitem{Liu-Ma2018}
C.-H. Liu and J.~Ma.
\newblock Cycle lengths and minimum degree of graphs.
\newblock {\em J. Combin. Theory Ser. B}, 128:66--95, 2018.

\bibitem{connectivity-2021}
K.~S. Lyngsie and M.~Merker.
\newblock Cycle lengths modulo {$k$} in large 3-connected cubic graphs.
\newblock {\em Adv. Comb.}, pages Paper No. 3, 36, 2021.

\bibitem{average-degree-2016}
J.~Ma.
\newblock Cycles with consecutive odd lengths.
\newblock {\em European J. Combin.}, 52:74--78, 2016.

\bibitem{Thomassen-induced-cycles}
C.~Thomassen and B.~Toft.
\newblock Non-separating induced cycles in graphs.
\newblock {\em J. Combin. Theory Ser. B}, 31(2):199--224, 1981.

\end{thebibliography}
\end{document}